\documentclass[11pt,a4paper]{amsart}
\usepackage{amscd,amssymb,amsopn,amsmath,amsthm,graphics,amsfonts,enumerate,
verbatim,calc}
\usepackage{fancyhdr}
\headsep=1cm \numberwithin{equation}{section}
\newtheorem{theorem}{Theorem}[section]
\newtheorem{proposition}[theorem]{Proposition}

\numberwithin{equation}{section}

\begin{document}
\thispagestyle{plain}
\begin{center}
{\large To appear in \sc Romanian Journal of Mathematics and Computer Science}

{\small available online at https:/$\!$/rjm-cs.utcb.ro }

{\footnotesize An Invited Address at RIGA 2025} \\
{\footnotesize Delivered in Bucharest on May 24, 2025} \\
\bigskip
{\small {\it Proceedings of the International Conference Riemannian Geometry and Applications} ‐ RIGA
2025} \\
{\small Bucharest, Romania, May 23-25, 2025}
\end{center}
\vspace{5cc}
\title{Extrinsic Characterisations of Immersions}
\author{Bogdan D. Suceavă}
\address{Department of Mathematics, California State University Fullerton, Fullerton, CA 92831-6850.}
\email{ bsuceava@fullerton.edu}

\maketitle

\begin{abstract}
In 1763, Leonhard Euler wrote that one cannot define a good curvature measure for surfaces:  ``la question sur la courbure des surfaces n'est pas susceptible d'une r\'eponse simple, mais elle exige \`a la fois une infinit\'e de d\'eterminations.'' The quest for the right measure for curvature was settled by C.F. Gauss in 1825, and Sophie Germain introduced the mean curvature in 1831 (her memoir written in 1816 included also the average of principal curvature as a shape invariant). We outline the history of the idea of deformation of space, which lead to the concept of curvature invariants, as we understand them today, including contributions of E. Bacaloglu and F. Casorati, among others. We pursue the following question: what is the best way to quantify the deformation of space? This important question could be viewed in a new paradigm after 1956, when John F. Nash, Jr. proved that a Riemannian manifold can be immersed isometrically into an Euclidean ambient space of dimension sufficiently large. This important theorem allowed to view the representation of space from its exterior, from an outside perspective. In 1968, S.-S. Chern pointed out that a key technical element in applying Nash's Theorem effectively is finding useful relationships between intrinsic and extrinsic quantities characterising immersions.  And such relations seem to be rather few, at least few enough to present us with a technical challenge in applying Nash's Theorem. One technical difficulty is presented by the passing through the narrow gateway provided by Gauss' equation, and that's why in order to obtain some new results it might be useful to include additional natural conditions. A turning point in the history of the question we pursue was an enlightening paper written by B.-Y. Chen in 1993, which paved the way for a deeper understanding of the meaning of the Riemannian inequalities between intrinsic and extrinsic quantities. Our present discussion invites a reflection on whether we could hope to characterise submanifolds by using mainly extrinsic quantities. After looking at several recent examples of such results, we conclude our paper with the following {\bf Problem.} {\it Are there any other extrinsic relations that determine the topology or the geometry of an embedded geometric object? How do we define them? How much insight do they provide, when we look at the geometric object ``from the outside"?}  \vspace{0.2 cm} \\

\end{abstract}

\vspace{.5cm}
{\it Dedication:} This work is dedicated to the memory of Professor Solomon Marcus (1925--2016), who pointed out that the history of curvature is incomplete without including E. Bacaloglu's work. In 2025 the UNESCO Commission in Romania celebrates the Solomon Marcus centenary.
\vspace{.5cm}

\section{The Question of Space and its Deformations}

The present paper represents the written version of one of the four Invited Addresses at the Riemannian Geometry and Its Applications (RIGA) 2025 conference. The talk was presented on May 24, 2025. While the author is profoundly grateful to the organisers for their invitation, this work aims to reach out towards an interesting idea much indebted to the advances ushered in by John F. Nash, Jr.'s Embedding Theorem and its consequences (for its importance see e.g. \cite{C2011}). We felt this was the best opportunity to present a more thorough overview of our research direction, starting with the profound philosophical motivation of the theme. We have been much encouraged towards such a perspective by an inspiring paper, Leopold Verstraelen's work \cite{V2020}, an author whose views shaped much of our current interests. Our present perspective was preceded by yet another Invited Addres, delivered in the fall of 2022, for a national conference in philosophy in Romania, which can be read in the work \cite{S2023}. 

Quite unexpectedly, the definition of curvature appears for the first time at the middle of the 14th century, in the works of Nicole Oresme \cite{O1351}, see also \cite{SS2015} \footnote{The author of the present paper found out about Oresme's works from Bang-Yen Chen's work \cite{C2000}, and also from several private conversations about the origin of the concept of curvature that took place at Michigan State University about the time the survey \cite{C2000} was published. The atmosphere surrounding the Differential Geometry Seminar at Michigan State University in those years is recalled in the book \cite{S2018}, and the careful reading of \cite{O1351} lead among other works to \cite{SS2015,SV2023}.}.

It is quite remarkable to point out that Oresme’s reasons to introduce his concept of {\it curvitas} (as it makes perfect sense within the context of his {\it doctrine}) are not related to the mathematical motivations that later authors pursued to investigate deformations of shape, and together with it the modern concept of curvature. It is an idea born in Western Europe in the intellectual context of the Rediscovery of Aristotle \cite{R2004}, a fertile context different from the scientific revolution from the 17th century \cite{SV2023}, when the later mathematical discoveries took place. The paradigms were so different that it really requires a different type of mathematical reading.

The assertion we propose is that Oresme introduced the idea of investigating shape (a thought encoded in his works under the far-enveloping word “configurations”) at the same time with introducing the idea of curvature as deformation of shape, which makes the so-called (by its author) {\it doctrine of configurations} a striking contribution to the medieval thinking, analogue in this aspect with Riemann’s fundamental contribution from 1854 \cite{R1854}. Namely, Oresme considers: (1) a representation of the idea of space; together with (2) its deformation. This is how there is an analogy between Oresme's {\it doctrine} and Riemann's seminal construction of the concept of manifold, together with its  sectional curvature \cite{R1854}.
 
 If there is any other source where the concept of curvature is discussed prior to 1351, it did not reach us and it is not known to us. It would actually be a very interesting question for historians if anyone would ever discover any contribution about the idea of deformation of space prior to Oresme's thoughts.

We cannot separate any discussion of Oresme’s doctrine of configurations from the problem of representing space. The origin of the theme, however, has in Western philosophy older roots than the Aristotle’s tradition, where Oresme found his inspiration (as far as his testimonial records). Bertrand Russell reflects in \cite{R1945}, p.146 et al., on the following key fragment, from Plato’s {\it Timaeus}: ``There is one kind of being which is always the same, uncreated and indestructible, never receiving anything into itself from without, not itself going out to any other, but invisible and imperceptible by any sense, and of which the contemplation is granted to intelligence only. And there is another nature of the same name with it, and like to it, perceived by sense, created, always in motion, becoming in place and again vanishing out of place, which is apprehended by opinion and sense. And there is a third nature, which is  space, and is eternal, and admits not of destruction and provides a home for all created things,  and is apprehended without the help of sense, by a kind of spurious reason, and is hardly real; which we beholding as in a dream, say of all existence that it must of necessity be in some place and occupy a space, but that what is neither in heaven nor in earth has no existence."

We believe that it might useful for anyone interested in geometry today to start any examination of geometric ideas with this highly inspiring starting point.  About this passage, Bertrand Russell writes \cite{R1945}:  ``This is a very difficult passage, which I do not pretend to understand at all fully. The theory expressed must, I think, have arisen from reflection on geometry, which appeared to be a matter of pure reason, like arithmetic, and yet had to do with space, which was an aspect of the sensible world. In general it is fanciful to find analogies with later philosophers, but I cannot help thinking that Kant must have liked this view of space, as one having an affinity with his own."

From the philosophical standpoint (see e.g. \cite{S2023}), therefore, the question we discuss therefore could be phrased as: “what is space?”  We argue that any discussion on the history of this theme should include Nicole Oresme’s contributions, even if we are hereby proposing an author who did not use a precise computational tool to determine curvature for planar curves, space curves, or for surfaces, as later authors did, either by using the tools of calculus, or later on vectorial or tensorial calculus. Nevertheless, for this particular topic, Oresme’s contribution was way ahead of his time \cite{S1969}.

\section{From Leonhard Euler to C.F. Gauss and Sophie Germain}

One generation before Gauss, the question of curvature of a surface at a given point was far from being settled. To better understand how the concept of curvature of surfaces was viewed before Gauss, we recall that Leonhard Euler wrote in 1763 that we cannot define the curvature for surfaces. He writes:

\begin{quote}
Pour connoitre la courbure des lignes courbes, la détermination du rayon osculateur en fournit la plus juste mesure, en nous présentant pour chaque point de la courbe un cercle, dont la courbure est précisément la même. Mais, quand on demande la courbure d'une surface, la question est fort équivoque, et point du tout susceptible d'une réponse absolue, comme dans le cas précédent.  
\end{quote}

The geometric situation that eluded Euler was the case when the tangent plane at a point to a surface is intersecting the surface, that is, if we were to use a more contemporary language, the case when the surface is hyperbolic at that point. This happens, e.g., at every point of the catenoid. What Euler hoped to see materialised was the idea of approximation of the surface with a sphere, an extension of Newton’s idea of osculating circle approximating locally the curve. However, it was not this idea that produced the definition of the curvature at a point of a surface, as we can see by reading Gauss’  {\it Disquisitiones} but a combination of ideas including small infinitesimals and limiting processes. It was Gauss’ vision that generated the curvature of surfaces, and later on it was Riemann's effort who extended this concept to $n$-dimensional representation of space.

In section V of  {\it Disquisitiones}, Gauss introduces the curvature for surfaces. 
For a smooth surface $S$ lying in $\mathbb{R}^3,$ and an arbitrary point $P \in S,$  consider $N_P$ the normal to the surface at $P.$  
Consider the family of all planes passing through $P$ that contain the line through $P$ with the same direction as $N_P.$  
These planes yield a family of curves on $S$ called {\it normal sections}. 
Consider now the curvature $\kappa(P)$ of the normal sections, viewed just as planar curves. Then $\kappa(P)$ has a maximum, denoted $\kappa_1$, and a minimum, denoted $\kappa_2.$    The curvatures $\kappa_1$ and $\kappa_2$ are called the {\it principal curvatures}.  The Gaussian curvature \cite{G} is defined as $K(P)= \kappa_1(P) \cdot \kappa_2(P),$  and Sophie Germain's mean curvature is defined \cite{Germain1831} by the arithmetic mean $H(P)= \frac{1}{2}\left[  \kappa_1(P) + \kappa_2(P) \right].$  C.F. Gauss' {\it Theorema Egregium}, a result of massive importance in geometry's development, was later rephrased by many authors that the curvature $K$ is of intrinsic nature, which means that it depends only of the geometry of the surface, and it is derived based on Gauss's calculations at the end of section 11 in  {\it Disquisitiones}.

It is quite interesting to read that Sophie Germain appreciated Euler's contribution as fundamental. She writes:

\begin{quote}
Si, par rapport aux surfaces, on avoit besoin de connoître la mesure de la courbure, on trouveroit peu de secours dans les \'ecrits des g\'eom\`etres qui se sont occup\'es des diverses questions dont se compose la g\'eom\'etrie descriptive, et l'on seroit forc\'e de remonter \`a des travaux plus anciens.

La m\'emoire d'Euler, intitul\'e: {\it Recherches sur la courbure des surfaces} contient, en effet, tout ce que l'on sait d'important \`a cet \'egard; et, en lisant ce beau m\'emoire, on apper\c{c}oit bient\^ot que l'illustre auteur y a d\'epos\'e le germe des recherches, qui peuvent faire disparoître les difficult\'es qu'il a pris soin de signaler.

\end{quote}
At p.4 in \cite{Germain1831}, S. Germain mentions that the mean curvature was included in her memoire submitted to the French Academy in 1816.

\section{Nicolas Renard and Emanoil Bacaloglu}

A noteworthy contribution in investigating the concept of curvature of surfaces took place on August 18, 1856, when at the Facult\'e des sciences de Paris Nicolas Renard defended his doctoral thesis, which was straightforwardly titled  {\it Courbure des surfaces} \cite{R1856}, with a committee chaired by Louis Lef\'ebure de Fourcy, and having as members Gabriel Lam\'e and  Augustin-Louis Cauchy (the actual defense took place just a few months before Cauchy passed away). Nicolas Renard investigation focused on the so-called geometric locus of the centers of mean curvature. He investigated the relationship between a given smooth surface and the surface which is this locus, and discovered that the two surfaces are reciprocal, in the sense introduced by him.

Renard's thesis is one of the references in Emanoil Bacaloglu's short article \cite{B1859}, in which he proposes a new curvature invariant. After a reasoning relying on a limiting process which involves Gauss' map, Bacaloglu argues that it would be natural to consider as curvature the quantity $H^2 + \frac{1}{8}L^2$, where $L$ is the difference between the principal curvatures at a point of a smooth surface. Bacaloglu's curvature invariant is an extrinsic concept. His vision regards the surface from the outside, in the paradigm proposed by Sophie Germain, although his reasoning starts from Gauss' map. What is truly important in this early development is the very connection between limiting processes and a quantity that {\it is}, in every sense of the word, representing the deformation of space. And it's not just that, Bacaloglu investigates, just for elliptical surfaces, for situations when $K>0,$ yet another limiting process that ends up with another quantity representing curvature, namely the invariant $K^{3/2}/ H.$ 

It was a very interesting early development in the geometry of smooth surfaces, showing that one indeed can consider rightfully other quantities other than $K$ or $H$, as long as the quest is motivated by a profound geometric substance.

\section{Felice Casorati's Curvature Invariant. What Can Be Pursued Today?}

In a recent very interesting paper \cite{RP2024}, Arrigo Pisati and Riccardo Rosso discuss the genesis of Felice Casorati's paper from 1890 \cite{C1890}, where he proposed his curvature invariant. Pisati and Rosso self-describe their effort: ``In his last paper, Felice Casorati proposed an extrinsic measure of curvature of a surface that vanishes only on planes. His proposal triggered strong reactions in the mathematical community that we examine here, also by use of unpublished material." Actually, Casorati's inquiries towards a better understanding of the idea of deformation of space go back to 1867, see \cite{C1867}. We feel it is important to understand well Casorati's argument and his efforts, and a lecture of his original work is not lacked of interest today. It is truly the spirit of that time  (the decade of 1890s) that Casorati's argument is fundamentally a limiting process motivated by geometric intuition. Casorati proposes as curvature of a surface at a point $p$ the quantity $C(p) = \frac{1}{2}[ \kappa_1(p)^2 + \kappa_2(p)^2],$ and in his paper \cite{C1890} we read plenty of reasons to do so. Casorati, as well as Bacaloglu back in 1859, has as starting point a limiting process inspired by Gauss' profound reflection dating back to 1825 \cite{G}. Pisati and Rosso are right to point out that recently there is a new revived interest in investigations in which the Casorati curvature is regarded as a geometric quantity representing in the most legit way the idea of curvature, see e.g. some recent works in this direction \cite{ALVY, BS2019, C2021,CDV2021, D2014, D2007, D2008,LLV2020,LV2015,L2017,SV2019,V2018}, to recall here just a few.

For the effectiveness of our presentation, we feel it would be best if we pursue our analysis in the context of smooth hypersurfaces embedded in a real Euclidean ambient. The reason we do so is that we would like to connect this historical heritage to the contemporary theory of submanifolds. To recall a few concepts in the geometry of differential hypersurfaces or, to be more precise, smooth hypersurfaces \cite{BS2018,S2015,SV2019}, let $\sigma : U \subset \mathbb{R}^n \rightarrow \mathbb{R}^{n+1}$ be a smooth hypersurface given by the smooth map $\sigma.$ Let $p$ be a point on the hypersurface. Denote  $\sigma_k(p) = \frac{\partial \sigma}{\partial x_k},$ for all $k$ from 1 to $n.$ Consider $\{ \sigma_1(p), \sigma_2(p), ..., \sigma_n(p) , N(p) \},$ the Gauss frame of the hypersurface, where $N$ denotes the normal vector field. We denote by $g_{ij}(p)$ the coefficients of the first fundamental form and by $h_{ij}(p)$ the coefficients of the second fundamental form. Then 
$$ g_{ij}(p) = \langle \sigma_i(p), \sigma_j(p) \rangle, \  \  \  \ \  h_{ij}(p) = \langle N(p), \sigma_{ij} (p)\rangle.$$

The Weingarten map $L_p = - dN_p \circ d\sigma_p^{-1} : T_{\sigma(p)} \sigma \rightarrow T_{\sigma(p)} \sigma $ is linear. Denote by $( h_j^i(p))_{1\leq i,j \leq n}$ the matrix associated to Weingarten's map, that is:
$$L_p( \sigma_i(p) = h^k_i(p) \sigma_k(p),$$
where the repeated index and upper script above indicates Einstein's summation convention. Weingerten's operator is self-adjoint, which implies that the roots of the algebraic equation
$$ \det ( h_j^i(p) - \lambda (p) \delta_j^i ) =0$$
are real. The eigenvalues of Weingarten's linear map are called principal curvatures of the hypersurface. They are the roots $k_1(p), k_2(p), ... , k_n(p)$ of this algebraic equation.
The mean curvature at the point $p$ is $$H(p) = \frac{1}{n} [k_1 (p) + ... +k_n(p)],$$ and the Gauss-Kronecker curvature is $$K(p) = k_1 (p) k_2(p) ... k_n(p).$$

If all the principal curvature of a smooth regular hypersurface are $\geq 0,$ then the hypersurface is convex.

One may define the Casorati curvature in the direction of the $k$-dimensional planar section $V_p \subset T_pM,$ or the Casorati curvature of order $k$, by $$C_k(e_1, e_2,...,e_k)= a_1^2+ ... + a_k^2,$$ for $k \leq n,$ where $a_1,...,a_k$ are the principal curvatures at point $p \in M.$ If we regard the matters in this way, Casorati curvatures are related to the geometry of submanifolds \cite{C1973,C1981,C2011}, and any information that relates curvature invariants to the topology of the submanifolds is of interest. Here is the idea we pursue. 

{\bf Question.} {\it Are there any prescribed conditions in terms of Casorati curvature that yield global conclusions about the geometry of a surface or of a hypersurface?}

To give an example of what we can obtain from this idea, we refer to the following.

\begin{proposition} \cite{BS2018}
\label{BS}
Let $\sigma: U \subset \mathbb{R}^{n} \rightarrow \mathbb{R}^{n+1}$ be a regular smooth hypersurface, $Im \ \sigma = M^n.$ Let $a_1, a_2,..., a_n$ be the principal curvatures at $p \in M.$ If all the Casorati curvatures of order 2 satisfy the inequality
$$\sqrt{(n-1) C_{n-1}(p)} \leq nH(p),$$
for every $p$ in $M,$ then the hypersurface must be convex.
\end{proposition}

This result is cited also in \cite{C2021}, where the derivation of this analytic fact is placed in a natural context of ideas.

For its proof, we need the following claim  \cite{BS2018}.
 If for $n \geq 3$ real numbers $a_1,a_2,...,a_n$ the following inequalities hold:
$$ \sqrt{(n-1) (a_1^2 +a_2^2+ ... + a^2_{n-1})} \leq a_1+a_2 +...+a_n,$$
$$ \sqrt{(n-1) (a_1^2 +a_2^2+...+a^2_{n-2}+ a^2_n)} \leq a_1+a_2 +...+a_n,$$
$$\ldots$$
$$ \sqrt{(n-1) (a_2^2 +a^2_3+...a_n^2)} \leq a_1+a_2 +...+a_n.$$
Then we must have $a_1 \geq 0, a_2 \geq 0, $ ... $a_n \geq 0.$

To prove this claim, we solve in the first inequality for $a_n$, which appears only in the right hand side term:
$$a_n \geq  \sqrt{(n-1) (a_1^2 +a_2^2+ ... + a^2_{n-1})} - (a_1+a_2 +...+a_{n-1}).$$
We need to prove this term is greater than or at least equal to zero.
Then we have:
$$a_n \geq  \sqrt{(n-1) (a_1^2 +a_2^2+ ... + a^2_{n-1})} - (a_1+a_2 +...+a_{n-1}).$$
We need to prove the right side of this inequality is greater than or at least equal to zero. This is just the Cauchy-Schwarz inequality applied to the numbers $a_1,a_2 ,...,a_{n-1},$ and then to $(n-1)$ copies of $1:$
$$ \sqrt{(1^2+1^2+...+1^2) (a_1^2 +a_2^2+ ... + a^2_{n-1})} \geq (1\cdot a_1+ 1\cdot a_2 +...+1\cdot a_{n-1}).$$
This proves that $a_n \geq 0.$ Similarly one may prove all the other inequalities $a_i \geq 0, $ $i = 1, ..., n-1.$
Equality holds when all the $n$ numbers $a_i$, $i=1,...,n$ take the same value. 
 Proposition \ref{BS} is a direct consequence of this algebraic argument. \qed

The important point is this: The inequality in the above Proposition  \ref{BS}   \cite{BS2018} is an all-extrinsic condition, which yields the convexity of the geometric object under investigation. One can even draw some analogies with the Second Derivative Test from our calculus experience. In this context, we would like to conclude this section with the following open thought: {\it Are there any more results of this kind awaiting for us to discover them?} It may be just a matter of us adopting the appropriate perspective, as we follow the path and the paradigm opened by Nash's Embedding Theorem.

\section{A Question of S.-S. Chern}

Consider a surface $S$ in $\mathbb{R}^3.$ Principal curvatures at 
$P \in S$ are
$k_1, k_2.$ We have seen that the 
Gaussian curvature is $K(p) = k_1 \cdot k_2,$ and Sophie Germain's 
mean curvature is $H(p) = \frac{ k_1 + k_2}{2}.$

Remark that $ (k_1 - k_2)^2 \geq 0,$ or, if we rewrite   $k_1^2 - 2k_1k_2 + k_2^2 \geq 0,$ we obtain
$k_1^2 + 2k_1k_2 + k_2^2 \geq 4 k_1k_2.$
In terms of curvature invariants, at any point $P \in S$:
$$(k_1+k_2)^2 \geq 4 K.$$ 
In conclusion: $H^2(P) \geq K(P).$ We read this relation as an inequality of type {\it extrinsic} $\geq$ {\it intrinsic.}

Thus, the inequality $H^2(P) \geq K(P)$ is an elementary relationship of type {\it extrinsic} $\geq$ {\it intrinsic}.

Why is such a relationship important for the geometry of the surface embedded in the three dimensional real space?
A surface $S \subset \mathbb{R}^3$ is minimal if at $\forall P\in S: $ $H(P) \equiv 0.$ Since $H^2(P) \geq K(P),$ we  obtained the following {\it obstruction to minimality:} If there exists a point $P \in S$ such that  $K(P) >0,$ then $H(P) \neq 0.$

A classical theorem in the geometry of surfaces states that  if $S$ is compact, then there exists $P \in S$ such that $K(P) >0.$ Based on the above obstruction, we have just proved that a compact surface in $\mathbb{R}^3$  can not be minimal.

After this elementary take, we turn our attention to a more general case.
There is more than the sectional curvature (introduced first by B. Riemann, 1854).
  For an orthonormal
basis $e_1,...,e_n$ of the tangent space $T_p M$ of a Riemannian
$n$-manifold $M$, the scalar curvature at $p$ is:
$$ scal(p) =  \sum _{i<j} sec (e_i \wedge e_j ).$$
If we denote by
$h$ the second fundamental form of an isometric immersion of
a Riemannian $n$-manifold $M^{n}$ into a Riemannian space $\bar
M^{n+m},$ the 
mean curvature vector field by $H = \frac{1}{n}\,\hbox{ trace}
\, h,$ 
the immersion is called {\bf minimal} if $H(p) = 0,$ $\forall p \in M.$

Consider now $M,$ a submanifold of dimension $n$ in a Riemannian ambient space $(\bar{M}, \bar{g})$
of dimension $n+m.$

Gauss formula is:
$$\bar{\nabla}_XY = \nabla_XY + h(X,Y),  \ \ \  \forall X,Y \in T_pM.$$  
Gauss equation is:
$$g(R(X,Y)Z, W)= \bar{g}( \bar{R}(X,Y)Z,W) + \bar{g}(h(X,W),h(Y,Z)) $$
$$- \bar{g}(h(X,Z), h(Y,W)).$$   

In particular, if $M$ is a submanifold in an Euclidean ambient space:
$$sec ( e_i \wedge e_j) = \bar{g} (h(e_i,e_i), h(e_j,e_j)) - |h(e_i,e_j)|^2.$$

Let $v=e_1$ be a unit vector at $P.$ Consider $e_2,...,e_n$ the completion to an orthonormal basis of $v.$ Then, by definition:
$$Ric (v,v) = \sum_{i=2}^n sec( v \wedge e_i)$$  
By using Gauss equation: $$Ric (e_1,e_1) 
=\sum_{i=2}^n sec( e_1 \wedge e_i)=\sum_{i=2}^n \bar{g}(h(e_1,e_1) , h(e_i,e_i) )- \sum_{i=2}^n |h(e_1, e_i)|^2=$$  
$$= \bar{g}\left( h(e_1,e_1) , \sum_{i=2}^nh(e_i,e_i)\right) - \sum_{i=2}^n |h(e_1, e_i)|^2=$$ 
$$= \bar{g}\left( h(e_1,e_1) , nH - h(e_1,e_1)\right) - \sum_{i=2}^n |h(e_1, e_i)|^2= \bar{g}\left( h(e_1,e_1) , nH \right) - \sum_{i=1}^n |h(e_1, e_i)|^2.$$

Thus, we have seen that if $M$ is a minimal  submanifold ($H=0$ at every point)  in an Euclidean space, by Gauss equation we have
$$Ric(X,X) = - \sum_{i=1}^n |h(X,e_i)|^2 \leq 0,$$
where $ \{ e_1, ..., e_n \}$ is an orthonormal local frame field on $M.$
This means that  the Ricci tensor
of a minimal submanifold $M$ of a Euclidean space is negative
semi-definite. This is what Shiing-Shen Chern (1911-2004) pointed out in  his 1968 monograph titled
{\it Minimal submanifolds in a Riemannian manifold} \cite{C1968}, precisely at p.13,
that $Ric > 0$ yields the only known Riemannian obstruction to minimal isometric
immersion in Euclidean space. As Bang-Yen Chen pointed out many times, based on S.-S. Chern's remark, it is important to understand the Riemannian obstructions to minimality. (Or to any other classes of immersions, as much as they can be controlled through the geometric structures.)

Related to Chern's remark from  \cite{C1968}, Bang-Yen Chen asked and thoroughly investigated in the 1990s the following {\bf Problem:} {\em When does a given Riemannian manifold $M$ admit (or
does not admit) a minimal immersion into a Euclidean space of arbitrary
dimension ?}

This is the context of ideas in which Bang-Yen Chen proved the following theorem.

\begin{theorem} \cite{C1993} Let $M^n$ be a submanifold in a space form of constant sectional curvature $c.$ Then
$$\inf (sec) \geq scal - \frac{n^2(n-2)}{2(n-1)} |H|^2 - \frac{(n+1)(n-2)}{2} c.$$
The equality case is completely determined by the form of the shape operators with respect to a suitable orthonormal frame fields. 
\end{theorem}

This geometric inequality sparked an important direction of investigation, much of which is surveyed and profoundly analysed in the monograph \cite{C2011}. Our reflections in the present paper are much inspired and indebted to this philosophy and vision.

\section{An Olympiad Problem}

Although the inclusion of a rather elementary (but nontrivial) example here might seem surprising, we believe this thought has simultaneously expository and mathematical value. At the 1984 International Mathematical Olympiad, Problem 1 asked to prove that for $a,b$ and $c$ nonnegative real numbers satisfying the constraint $a+b+c=1$ the following double inequality holds true: $$0 \leq ab+bc+ca- 2abc \leq \frac{7}{27}.$$ The problem's solution is available in many references, e.g. \cite{APS} or, for a more recent source, \cite{AS2021}, pp.11-12.

We are interested in it, because this double inequality admits a nice geometric interpretation, that could reveal interesting information on other inequalities as well. Further similar extensions or generalisations would also be interesting, that's why this problem deserves special consideration.

The 1984 IMO Problem recalled above is about three nonnegative real numbers $a,b,c.$ In order to provide the geometric interpretation we mentioned, we focus our attention on a hypersurface in the four-dimensional Euclidean space $\mathbb{R}^4.$ We regard the numbers $a,b,c$ as the principal curvatures at a point on this hypersurface, and we aim to decode the meaning in this double inequality.
Let $\sigma : U \subset \mathbb{R}^3 \rightarrow \mathbb{R}^{4}$ be the graph in the four-dimensional real space of a function; this is a hypersurface given by the smooth map $\sigma.$ Let $p$ be a point on the hypersurface.
 
There are quantities similar to  $\kappa_1$ and  $ \kappa_2$ from the geometry of surfaces; and they are the principal curvatures of the hypersurface, denoted  $\lambda_1, \lambda_2, \lambda_3$. They are introduced as the eigenvalues of the so-called Weingarten linear map. 
Similar to the geometry of surfaces, the curvature invariants in higher dimensions can also be described in terms of the principal curvatures. Let $p$ be a point of the hypersurface $M$ immersed into $\mathbb{R}^4$ endowed with the canonical metric. Let $e_1, e_2, e_3$ an orthonormal frame at $p,$ that diagonalizes the Weingarten operator. Then, as Bernhard Riemann saw it in 1854, the natural analogue of Gaussian curvature in the direction of the planar section generated by $e_i$ and $e_j$ is the product $ \lambda_i \lambda_j$, where $i,j \in \{ 1,2,3 \},$ $i \neq j.$ We write $sec (e_i \wedge e_j) = \lambda_i \lambda_j.$
The mean curvature at the point $p$ is $$H(p) = \frac{1}{3} [\lambda_1 (p) +\lambda_2(p) +\lambda_3(p)],$$ and the Gauss-Kronecker curvature is $$K(p) = \lambda_1 (p) \lambda_2(p) \lambda_3(p).$$ In Riemannian geometry, a third important curvature quantity is the scalar curvature (\cite{C2011}, p.19) denoted by $scal(p)$, which intuitively sums up all the sectional curvatures on all the faces of the trihedron formed by the tangent vectors in the Gauss frame:
$$scal(p) = \sec( \sigma_1 \wedge \sigma_2) + \sec (\sigma_2 \wedge \sigma_3) + \sec(\sigma_3 \wedge \sigma_1)= \lambda_1\lambda_2+ \lambda_3\lambda_1+ \lambda_2\lambda_3. $$
The last equality is due to the Gauss equation of the hypersurface $\sigma(U)$ in the ambient space $\mathbb{R}^4$ endowed with the Euclidean metric. 

In order to match the notation in the Problem recalled above, denote $\lambda_1=a,$ $ \lambda_2=b,$ $\lambda_3=c.$ The condition in the hypothesis that $a,b,c >0,$ is quite strong, and it means that the hypersurface is strictly convex (see \cite{KN1969}, p.40). Then the statement of the Problem could be reformulated as:  {\bf  If in the neighborhood of point $p$ on a strictly convex three dimensional hypersurface the mean curvature $H = \frac{1}{3}(a+b+c)$ is constant and equal to $\frac{1}{3}$, then we have
$$0 \leq scal(p)- 2K(p) \leq \frac{7}{27}.$$}

The matter could be pursued even further, and it has to do with the rephrased assertion $$scal(p) \leq 2K(p) + \frac{7}{27}.$$ This fact must be appreciated under its full limelight. The left-hand side quantity, $scal(p)$ is {\it intrinsic} (as sum of three sectional curvatures), while the right hand side, represented by the double of the Gauss-Kronecker quantity plus a constant, is for a three-dimensional hypersurface {\it extrinsic} (see \cite{KN1969}, p.33). The 1984 IMO problem describes nicely an inequality between an intrinsic and an extrinsic quantity, exactly the kind of statement that made Gauss reflect on the nature of a geometric quantity describing shape. We did use Gauss's equation in a subtle way, when we said that the product $ \lambda_i \lambda_j$ represents sectional curvature.

This fundamental example allows us to hereby introduce a major development in the theory. In 1956 John F. Nash, Jr. proved \cite{F1956} that a Riemannian manifold (the hypersurface is a particular case) can be immersed isometrically into an Euclidean ambient space of dimension sufficiently large. We recalled that over a decade later, 
S.-S. Chern pointed out (see e.g. p.13 in \cite{C1968}) that a key technical element in applying Nash's Theorem effectively is finding useful relationships between intrinsic and extrinsic quantities characterizing immersions. And such relations are few. The fact that for constant mean curvature convex hypersurfaces of dimension 3 there exists such a relationship, namely $scal(p) \leq 2K(p) + \frac{7}{27},$ is no little thing. It would be interesting to find out more such relations, if they exist. Anyways, an additional geometric assumption provided a fundamental inequality. 

It might seem surprising that some elementary inequalities could provide an insight into the problem of the best possible immersion of a space into another ambient space.  As about the hypersurfaces with constant mean curvature, much could be said, see e.g. Chapter 4 in \cite{C1973}, for a thorough coverage of the idea. The original condition $a+b+c=1$ admits an interesting geometric interpretation as well, which makes the matter even more natural to be investigated, from a geometric standpoint.

\section{B.-Y. Chen's Curvature Invariants}

As we outlined before, in his  visionary paper \cite{C1993}  B.-Y. Chen proved (more precisely in Lemma 3.2), that for a submanifold $M^n$ in a space form $R^{n+m}(c)$ of constant sectional curvature $c$ the scalar curvature satisfies at a point the fundamental inequality
\begin{equation}
\label{delta}
scal - \inf (sec) \leq \frac{n^2(n-2)}{2(n-1)} |H|^2 + \frac{(n+1)(n-2)}{2} c,
\end{equation}
where $|H|$ represents the magnitude of the mean curvature vector, and $\inf(sec)$ represents the infimum of all the scalar curvature taken over all 2-planes at that respective point.  Recall that 
for any orthonormal
basis $e_1,...,e_n$ of the tangent space $T_p M$ in a Riemannian manifold $M^n$, the scalar
curvature is defined to be $ scal(p) =  \sum _{i<j} sec
(e_i \wedge e_j ).$  Chen's inequality recalled above is important because it illustrates the kind of relationships it would be interesting to obtain: between intrinsic geometric quantities, by one hand (the terms in the left), and extrinsic geometric quantities by the other (in the right). This is the kind of relations we are interested in finding out. The comprehensive monograph \cite{C2011} describes a series of 
 relationships recently discovered (in the last three decades) involving intrinsic and extrinsic curvature invariants. And that is important, it represents exploiting fully the pathways opened by John Nash's Embedding Theorem \cite{F1956}.
 
   For an orthonormal
basis $e_1,...,e_n$ of the tangent space $T_p M$ of a Riemannian
$n$-manifold $M$, the scalar curvature at $p$ is:
$$ scal(p) =  \sum _{i<j} sec (e_i \wedge e_j )$$

For any $r$-dimensional subspace of $T_p M$ denoted by
$L$, with orthonormal basis $e_1,...,e_r,$ one  defines
$$ scal (L) = \sum _{1 \leq i < j \leq r} sec (e_i \wedge e_j).$$

  B. - Y.  Chen considered the finite set $S(n)$ of
$k$-tuples $(n_1,...,n_k)$ with $k\geq 0$ which
satisfy the conditions:    $n_1 < n,$  $n_i\geq 2$   
and $n_1+...+n_k\leq n.$   
For each $(n_1,...,n_k) \in S(n)$ B. - Y.  Chen
introduced  the following Riemannian invariants:    
$$
      \delta (n_1,...,n_k)(p) = scal (p) - \inf \{
scal (L_1)+...+ scal (L_k) \}(p),
$$
where infimum is taken for all possible choices of orthogonal subspaces
$L_1,...,L_k,$ satisfying $n_j = $dim $L_j$, $(j=1,...,k).$

These invariants are now known as Chen invariants.  
When $k=0$,  the Chen invariant is nothing but the
scalar curvature.

With the following notation,
$$ c(n_1,...,n_k) = \frac{n^2(n+k-1-\sum n_j)}{2(n+k-\sum n_j)},$$
$$ b(n_1,...,n_k)=\frac{1}{2} \left\{ (n(n-1) - \sum _{j=1}^k n_j(n_j -
1)\right\},$$
we can recall the following.

\begin{theorem} [B.-Y. Chen's fundamental inequalities, Thm.13.3/pg.262 in \cite{C2011}]
 \label{basic} For any  $n$-dimensional
submanifold $M$ of an arbitrary Riemannian space 
$\bar{M}^{n+m}$ and for any  $k$-tuple
$(n_1,...,n_k) \in S(n)$, we have:
$$
     \delta(n_1,...,n_k) \leq c(n_1,...,n_k)|H|^2 + b(n_1,...,n_k)
\max \bar{sec}(p),
$$
where $\max \bar{sec}(p)$ denotes the maximum of the sectional curvature function of $\bar{M}^{n+m}$ restricted to
2-plane sections of $T_pM.$
Equality case completely determined by some conditions on the second fundamental form.
\end{theorem}

This important theorem was much cited and received a lot of attention from many geometers. One of its first consequences is the following.

\begin{theorem}[Obstruction to minimality, B.-Y. Chen, 2000 \cite{C2011}]
 \label{obs}
Let M be a Riemannian n-manifold. If there exists a k-tuple $(n_1,...,n_k) \in
S(n)$ and a point $p\in M$ such that
\begin{equation}
\delta (n_1,...,n_k)(p) > \frac{1}{2} \{ n(n-1) - \sum n_j(n_j - 1) \} \varepsilon,
\end{equation}
then M admits no minimal isometric immersion into any
Riem.space form $R^m(\varepsilon )$ with arbitrary codimension.

In particular, if
$\delta(n_1,...,n_k)(p) > 0$ at a point for some k-tuple
$(n_1,...,n_k) \in S(n),$ then M admits no minimal
isometric immersion into any Euclidean space for any codimension.
\end{theorem}

We presented earlier how the inequality $H^2(P) \geq K(P)$ in the geometry of surfaces yields a natural obstruction to minimality. This last result by B.-Y. Chen we recalled above represents a more advanced conquest of this idea, and a meaningful step forward in our better understanding the phenomenon.

\section{Extrinsic vs. Extrinsic Relations. An Open Ending}

Recently, the focus on extrinsic geometric quantities lead to the following fact. We are getting closer to the point we are making with the present paper.

\begin{theorem} \cite{S2023b}
\label{FCSB}
Let $p_0$ be an umbilic on a smooth hypersurface in $\mathbb{R}^{n+1},$ endowed with a second fundamental form denoted by $h$, whose averaged trace yields the mean curvature $H$. Let $s(L)$ be the spread of the shape operator. Let $p$ be another point close enough to  $p_0$, that is not an umbilic. Then, the limit
\begin{equation}
\lim_{p\rightarrow p_0} \frac{s(L)}{\sqrt{||h||^2 - nH^2}} \in \left[ \frac{2}{\sqrt{n} }, \sqrt{2} \right],
\end{equation}
provided the limit exists.
\end{theorem}

This assertion truly describes a limiting process with extrinsic quantity, along the lines of the philosophical perspective opened up by the paradigm of extrinsic-based invariants whose origin can be traced back to Sophie Germain (1816), Emanoil Bacaloglu (1859) and Felice Casorati's reflections (1890), as we described before. It is exactly this kind of results we are pursuing, due to this interesting tension of ideas and intertwining of perspectives which could be followed all along the history of differential geometry.

If  $s(L)$ represents the spread of the shape operator, the work in \cite{S2023b}  leads to the following inequality between intrinsic (in the left hand side term) and extrinsic quantities (grouped in the right hand side terms), respectively.

\begin{theorem} \cite{S2023b}
 \label{Steaua} Let $M^n \subset R^{n+1}(c)$ be a smooth hypersurface in a space form endowed with constant sectional curvature $c.$ If $||h||$ is the norm of the second fundamental form and $s(L)$ is the spread of the shape operator, then between the intrinsic and the extrinsic quantities at one point of the hypersurface the following inequality holds:
\begin{equation}
scal - \inf (sec) \leq \frac{n(n-2)}{4(n-1)} \left[ 2||h||^2 - s^2(L)] \right]+  \frac{(n+1)(n-2)}{2} c.
\end{equation}
Equality holds at the umbilical points.
\end{theorem}

For a sketch of the proof, we remark that at any point of the smooth hypersurface, we have:
\begin{equation}
\label{fact1}
 \delta(2) -  \frac{(n+1)(n-2)}{2} c \leq \frac{n^2(n-2)}{2(n-1)} H^2.
 \end{equation}

On the other hand,  $s^2(L) \leq 2 ||h||^2 - 2nH^2,$ which yields right away
\begin{equation}
\label{fact2}
H^2 \leq \frac{ 2||h||^2 - s^2(L)}{2n}.
\end{equation}
By combining the relations we obtain the claimed inequality.  Note that the umbilical points satisfy the equality case as described in Lemma 3.2 from B.-Y.Chen's fundamental work \cite{C1993}. \qed

We have the following consequence.

\begin{proposition} \cite{S2023b}
\label{rho}
Let $M^n$ be a hypersurface in a Riemannian
$(n+1)$-manifold $\bar M^{n+1}$. Then at every point
$p \in M$ the following inequality holds:
\begin{equation}
\rho(p) \leq \frac{ 2||h||^2 - s^2(L)}{2n}  + \frac{2}{n(n-1)}  \sum_{i<j} \overline {sec} (e_i\wedge
e_j),
\end{equation}
where $\rho$ is the normalized scalar curvature of $M$ at $p,$ $H$ is the
mean curvature at $p$, and
$\overline {sec} (e_i \wedge e_j)$ represents the sectional
curvature on the plane
generated by vectors $e_i$ and $e_j$ tangent to the ambient space
$\bar M$.
The equality holds at $p$ if and
only if $p$ is an umbilical point.
\end{proposition}

As an interesting point to make here, recall that in \cite{S2015} a class of purely extrinsic invariants have been introduced. The construction is the following.

\begin{theorem}
\label{second-dim3}
 Let $M^3 \subset \mathbb{R}^4$ be a smooth hypersurface and $k_1, k_2, k_3$ be its principal curvatures in the ambient space $\mathbb{R}^4$ endowed with the canonical metric. Let $p\in M$ be an arbitrary point. Suppose that at $p$ none of the principal curvatures vanish. Denote by $K$ the Gauss-Kronecker curvature. Introduce the
 absolute mean sectional curvatures defined by
 $$\bar{H}_{ij}(p) = \frac{|k_i(p)| + |k_j(p)|}{2}.$$ 
Then the following inequalities hold true at every point $p \in M:$
\begin{equation}
\label{Nesbitt}
3 \leq  \frac{|k_1|}{\bar{H}_{23} }+ \frac{|k_2|}{\bar{H}_{13}}+ \frac{|k_3|}{\bar{H}_{12}}   = B_1^1(p)
\end{equation}

\begin{equation}
\label{doi}
2\sqrt{2} < \sqrt{\frac{|k_1|}{\bar{H}_{23} }}+\sqrt{ \frac{|k_2|}{\bar{H}_{13}}}+ \sqrt{\frac{|k_3|}{\bar{H}_{12} }}  = B_{0.5}^{0.5} (p) 
\end{equation}

\end{theorem}
The first of these inequalities is nothing else but the classical Nesbitt's inequality, this time turned into a geometric fact. We have therefore the {\it geometric interpretation of Nesbitt's inequality}, as earlier in this paper we discussed a geometric interpretation of Cauchy-Schwarz inequality in terms of Casorati curvature. We thought this result deserves to be recalled here because it shows 
 that there are examples of extrinsic curvature invariants bounded below by a constant; their construction is provided in \cite{S2015}. 

\vspace{.2cm}

We would like to conclude our present work with the statement of the following.

\vspace{.2cm}

{\bf Problem.} {\it Are there any other extrinsic relations that determine the topology or the geometry of an embedded geometric object? How do we define them? How much insight do they provide, when we look at the geometric object ``from the outside"?}

\vspace{.2cm}

This way of asking the question is philosophically indebted to John Forbes Nash's Embedding Theorem, exactly as described by Bang-Yen Chen in \cite{C2011}. By thinking along these lines we are pursuing a classical direction of inquiry in differential geometry, which might have further ramifications, as the interest in a variety of geometric objects might diversity in the future, and whose noble history was hereby outlined \cite{B1859,C1867,C1890,G,Germain1831,E1767,O1351,R1854}.

\end{document}